\documentclass[eqnumis, pdf]{hjms}
\usepackage{amssymb}
\usepackage{lscape}
\usepackage{enumitem}
\usepackage{amsfonts}
\usepackage{amsmath}
\usepackage{footnote}
\usepackage{multicol}
\usepackage{booktabs}
\usepackage[flushleft]{threeparttable}
\usepackage[font=small,labelfont=bf]{caption}

\begin{document}
\hinfo{XX}{x}{XXXX}{1}{\lpage}{10.15672/hujms.xx}{Article Type}
%

\markboth{A. M. Mathai, M. Kumar}{Variable ASP for expo. dist. using Bayesian estimate under Type I hybrid censoring}

\title{Optimal variable acceptance sampling plan for exponential distribution using Bayesian estimate under Type I hybrid censoring}

\author{A. M. Mathai$^{1}$, M. Kumar\coraut$^2$}

\address{$^{1,2}$Department of Mathematics, National Institute of Technology, Calicut, Pin Code - 673601, India.}
\emails{ashlyn\_p190073ma@nitc.ac.in (A. M. Mathai), mahesh@nitc.ac.in (M. Kumar)}
\maketitle
\begin{abstract}
	In this study, variable acceptance sampling plans under Type I hybrid censoring is designed for a lot of independent and identical units with exponential lifetimes using Bayesian estimate of the parameter $\vartheta$. This approach is new from the conventional methods in acceptance sampling plan which relay on maximum likelihood estimate and minimising of Bayes risk. Bayesian estimate is obtained using squared error loss and Linex loss functions. Optimisation problem is solved for minimising the testing cost under each methods and optimal values of the plan parameters $n, t_1$ and $t_2$ are calculated. The proposed plans are illustrated using various examples and a real life case study is also conducted. Expected testing cost of the sampling plan obtained using squared error loss function is much lower than the cost of existing plans using maximum likelihood estimate.
\end{abstract}

\hinfoo{DD.MM.YYYY}{DD.MM.YYYY} 


\section{Introduction}\label{sec1}
In the world of manufacturing and quality control, ensuring the consistent quality of products is of paramount importance. However, inspecting each and every item produced can be impractical, time-consuming, and costly. This is where acceptance sampling plans come into play. An acceptance sampling plan (ASP) is a statistical technique used to make informed decisions about accepting or rejecting a batch or lot of items based on the inspection of a sample.

 Understanding and implementing a well-designed acceptance sampling plan is crucial for manufacturers seeking to maintain consistent product quality, optimize resources, and meet customer expectations. By employing statistically sound sampling techniques, organizations can make reliable decisions regarding batch acceptance, leading to improved efficiency, reduced costs, and enhanced customer satisfaction.

  Censoring schemes play a crucial role in ASP, and are designed to optimize the efficiency and effectiveness of the sampling process. Type-I and Type-II censoring schemes are widely used in statistical analysis, particularly in reliability testing and acceptance sampling. In a Type-I censoring scheme, the time of experiment is predetermined, but the total number of failures occurred becomes a random variable. Conversely, in a Type-II censoring scheme, the testing time is a random variable since the number of observed failures is fixed. When a combination of both Type-I and Type-II censoring schemes is employed, it is referred to as a hybrid censoring scheme (HCS).

Here, we consider an experiment where a sample of $n$ units are subjected to testing. Also, assume that the sample units' lifespans are independent random variables and identically follows exponential distribution. The sample units ordered lifetimes are represented as $X_{1,n}, X_{2,n}, ..., X_{n,n}$. The test is concluded either when a predetermined number, $\gamma$ (where $\gamma < n$), of the $n$ items have failed, or when a predetermined time, $\mathcal{T}$, has been reached. In other words, the life-test ends at a random time $\mathcal{T}^*$, which is the minimum of either the ordered lifetime ($X_{\gamma,n}$) or the predetermined time $\mathcal{T}$.  Moreover, a widely accepted assumption is that the failed units in the experiment are not substituted or replaced.

The Type-I HCS, proposed by Epstein \cite{epstein1954truncated}, has found extensive application in reliability acceptance testing, such as in MIL-STD-781-C (1977). 
Since its introduction, extensive research has been conducted on hybrid censoring and its various variations. Epstein initially introduced the Type-I HCS and discussed estimation methods for exponential distribution. Additionally, a two-sided confidence interval was put forth without a mathematical proof of how it was obtained.

Subsequent researchers, such as Fairbanks et al.\cite{fairbanks1982confidence}, Chen and Bhattacharyya \cite{chen1987exact}, Bartholomew \cite{bartholomew1963sampling}, and Barlow et al. \cite{barlow1968statistical}  made slight modifications to Epstein's proposition and worked on the derivation of the conditional moment generating function approach to derive the exact distribution of the conditional maximum likelihood estimator (MLE) of the mean life parameter . Childs et al. \cite{childs2003exact} obtained a simplified but equivalent form of the exact distribution of the MLE as derived by Chen and Bhattacharyya \cite{chen1987exact}. Draper and Guttman \cite{draper1987bayesian} delved into the Bayesian inference within the Type-I HCS and obtained Bayesian estimates and two-sided credible intervals using an inverted gamma prior.

In \cite{kumar2016design}, Kumar and Ramyamol, designed a cost efficient acceptance sampling plans for exponentially distributed lifetime data under Type I censoring based on MLE of the the mean life. Most of the works in acceptance sampling plans are centered on MLE of the parameter for deciding the acceptance and rejection of lot, even though Bayesian estimator is more reliable than MLE. 
 
Bayesian acceptance sampling plans have gained considerable attention due to their ability to incorporate prior information and make informed decisions based on observed data. However, existing Bayesian sampling plans focus on minimising Bayes risk, while the consideration of testing cost optimization remains limited. Lam \cite{yeh1994bayesian}, Lam and Choy \cite{yeh1995bayesian}, Lin et al. \cite{lin2002bayesian}, and Huang and Lin \cite{huang2002improved, huang2004bayesian} focused on conventional Type I, Type II, and random censored samples. Additionally, Chen et al. \cite{chen2004designing} examined Bayesian sampling plans for exponential distributions based on Type I hybrid censored samples.  

Recent contributions by Chen et al. \cite{chen2004designing}, Lin et al. \cite{lin2008exact} and Prajapati et al. \cite{prajapati2022optimal} have specifically addressed Bayesian sampling plans for exponential distributions by reducing Bayes risk utilizing Type I hybrid censored samples . But these plans are not based on Bayesian estimator of mean life parameter which is better than MLE (\cite{kundu2009estimating}). 

In this work, we examine an acceptance sampling plan for a batch of units where the failure time follows an exponential distribution with mean $\vartheta$. Our aim is to design an optimal Bayesian acceptance sampling plan under Type I HCS by using Bayesian estimator of the mean life parameter, $\vartheta$, as the decision function. Our plan stands out from all other existing Bayesian sampling plans by centering the decision of acceptance or rejection of a lot on the Bayesian estimator of the mean lifetime parameter. Also, here testing cost is minimized while maintaining stringent quality standards. This departure from traditional Bayes risk-oriented approaches allows us to focus explicitly on the optimization of testing costs, a critical consideration in practical industrial settings.

The rest of the paper is organised as follows: Section~\ref{sec2} provides a detailed information about Bayesian inference of exponential distributed life time data under Type I HCS. Bayesian estimator of mean life parameter, $\vartheta$, is considered using both squared error loss (SEL) function and Linex loss function. Section~\ref{sec3} presents the formulation and optimization framework of the proposed acceptance sampling plan. Here, the plan parameters ($n,t_1,t_2$) are obtained by minimising the expected testing cost (ETC) subject to required probability conditions using both the loss functions namely, SEL and Linex loss function. The distribution of the Bayesian estimator of $\vartheta$ under Type I HCS for both loss functions are derived using delta method. In Section~\ref{sec4}, we present the results of numerical computations of the proposed acceptance sampling plan and a comparison with existing work in \cite{kumar2016design}. A real life case study is done to illustrate the performance of the sampling plan presented here. Finally, Section~\ref{sec5} concludes the paper with a summary of the contributions and future research directions.  
\section{Bayesian Estimation for Type I hybrid censoring scheme}\label{sec2}
Consider a lot of units, all of which experience failures over time following an exponential distribution. The probability density function that characterizes the failure time of each unit is given by:
\begin{equation}
f(x,\vartheta)=\frac{1}{\vartheta}exp\left\{-\frac{-x}{\vartheta}\right\} \qquad x \geq 0, \vartheta > 0 
\label{eq1}
\end{equation}
 Here we conduct Type-I hybrid censored life testing, and we obtain one of the
following two forms of observations as our observed data:
\begin{equation}
    \begin{aligned}[b]
  \text{Data type I}: &\left\{X_{1,n} < X_{2,n} <\cdots<X_{\gamma,n}\right\} \quad \text{if} \quad \mathcal{T}^* = X_{\gamma},\\
\text{Data type II}: &\left\{X_{1,n} < X_{2,n}<\cdots<X_{\mathcal{D},n}\right\} \quad \text{if} \quad \mathcal{T}^* = \mathcal{T},       
    \end{aligned}
    \label{eq2}
\end{equation}
where $\mathcal{D}$ denotes the number of failures observed before time $\mathcal{T}^*$. Thus, $\mathcal{D}=\gamma$ when $\mathcal{T}^*=X_{\gamma}$.

In this section, the method of obtaining Bayesian estimator of the unknown parameter $\vartheta$ under the squared error loss and Linex loss functions and their corresponding distributions are discussed. Bayesian estimator of the parameter $\vartheta$ is used for the designing of acceptance sampling plan in the next sections and this idea is new in acceptance sampling plans. 

Under Type-I HCS, the likelihood function is obtained as 
\begin{equation}
L(x,\vartheta)=\frac{1}{\vartheta^\gamma}exp\left\{-\frac{1}{\vartheta}\sum\limits_{i=1}^{\mathcal{D}}x_{i,n}-\frac{(n-\gamma)}{\vartheta}\mathcal{T}^*\right\}, 
\label{eq3}
\end{equation} and 
  log-likelihood function as 
  \begin{equation*}
    L(\vartheta,x) = -\gamma \operatorname{log}\vartheta - \frac{1}{\vartheta} \sum\limits_{i=1}^{\mathcal{D}}x_{i,n}-\frac{(n-\gamma)}{\vartheta}\mathcal{T}^*,  
  \end{equation*}
  where $\mathcal{T}^*=\operatorname{Min}\{\mathcal{T},X_{\gamma,n}\}$ and
  $\mathcal{D}$ = Number of units failed before $\mathcal{T}^*$.
  
  Hence, by referring to the observed sample, the MLE of $\vartheta$  can be determined as follows (\cite{chen1987exact}):
  \begin{equation}
\hat{\vartheta}_{MLE}= \begin{cases}
      \frac{1}{\mathcal{D}}\left[\sum\limits_{i=1}^{\mathcal{D}}x_{i,n} + (n-\mathcal{D})\mathcal{T}\right] &  \text{if} \quad 1 \leq \mathcal{D} \leq  \gamma-1\\
      \frac{1}{\gamma}\left[\sum\limits_{i=1}^{\gamma} x_{i,n} + (n-\gamma)X_{\gamma,n}\right] &  \text{if} \quad  \mathcal{D} = \gamma\\
      n\mathcal{T} & \text{if} \quad \mathcal{D}=0.\end{cases}
      \label{eq4}
  \end{equation}

      It is evident that $\widehat{\vartheta}_{MLE}$ is a conditioned MLE of $\vartheta$, with the condition that at least one observed failure is present.

Here for Bayesian inference is done by using two kinds of loss functions, namely, squared error loss (SEL) and Linex loss function. SEL is a symmetric function and its symmetric nature is demonstrated as 
\begin{equation}
l_1(\eta,\widehat{\eta})=(\widehat{\eta}-\eta)^2 
\label{eq5}
\end{equation}
where $\widehat{\eta}$ is the estimate of the unknown parameter $\eta$.

The Bayesian estimation of any function $g=g(\vartheta)$ under the SEL in~(\ref{eq5}), is obtained by
\begin{equation}
\widehat{g}_{s}=E(g \mid Data)=\frac{\int_\vartheta g(\vartheta) L(x, \vartheta) \rho(\vartheta) d \vartheta}{\int_\vartheta L(x, \vartheta) \rho(\vartheta) d \vartheta}.
\label{eq6}
\end{equation}
where $\vartheta$ follows the prior distribution with PDF $\rho(\vartheta)$.\\

 Using Type-I HCS, Draper and Guttman \cite{draper1987bayesian} investigated the application of Bayesian inference to estimate the unknown parameter $\vartheta$ under squared error loss function. They assumed in the study that $\vartheta$ follows an inverted gamma prior distribution, with the probability density function (PDF) given below
      \begin{equation}
\rho(\vartheta)=\frac{a^b}{\Gamma(b)}\vartheta^{-(b+1)}e^{-a/\vartheta}, \quad \vartheta>0,\ a>0,\ \text{and} \ b>0. 
\label{eq7}
      \end{equation}
      The prior in~(\ref{eq7}) becomes non-informative prior when the hyper-parameters $a=b=0$. 

      The posterior density function of $\vartheta$, based on the knowledge from the observed data and the inverted gamma prior given by~(\ref{eq7}), becomes 
      \begin{equation}
          h(\vartheta \mid Data) = \frac{\left(\mathcal{D} \ \hat{\vartheta}_{MLE}+a\right)^{(\mathcal{D}+b)}}{\Gamma(\mathcal{D}+b)}\vartheta^{-\left(\mathcal{D}+b+1\right)} \operatorname{exp}\left\{\left(\mathcal{D} \ \hat{\vartheta}_{MLE}+a\right)/\vartheta\right\}
          \label{eq8}
      \end{equation}
      Considering that $\frac{\left(\mathcal{D} \ \hat{\vartheta}_{MLE}+a\right)}{\vartheta}$  is distributed as $\chi^2_{2\left(b+a\right)}/2$ a posteriori, the posterior mean serves as the Bayesian estimate of $\vartheta$ under the squared-error loss function. If $\mathcal{D}+b>1$, this estimator is simply produced as 
      \begin{equation}
\widehat{\vartheta_{s}} = \frac{\mathcal{D} \ \hat{\vartheta}_{MLE}+a}{\mathcal{D}+b-1}.
\label{eq9}
      \end{equation} When $a=b=0$, that is, when the prior becomes non-informative, it is evident that the MLE of $\vartheta$ provided by~(\ref{eq4}) agrees with the Bayesian estimator in~(\ref{eq9}).\\

The second is the Linex loss function, which is referred to as an asymmetric function and is represented by the following formula:
\begin{equation}
    l_2(\eta, \widehat{\eta}) \propto e^{c(\widehat{\eta}-\eta)}-c(\widehat{\eta}-\eta)-1, \quad c \neq 0.
    \label{eq10}
\end{equation}

The direction and magnitude of asymmetry are both determined by the parameter $c$'s value. Bayesian estimation tends to overestimate when $c<0$, but it favours underestimating when $c>0$. The Linex loss value resembles the squared error loss when $c$ approaches 0, leading to an almost symmetric behaviour. Applying the Linex loss function, we can derive the Bayesian estimation of $g(\vartheta)$ as
\begin{equation}
 \widehat{g}_l=-\frac{1}{c} \ln E\left(e^{-c g} \mid Data\right),
 \label{eq11}
\end{equation}
where
\begin{equation}
E\left(e^{-c g} \mid Data\right)=\frac{\int_\vartheta e^{-c g(\vartheta)} L(x, \vartheta) \rho(\vartheta) d \vartheta}{\int_\vartheta L(x, \vartheta) \rho(\vartheta) d \vartheta}. 
\label{eq12}
\end{equation}
But in most of the cases, the ratio of integrals in~(\ref{eq6}) and~(\ref{eq12}) is not possible to derive an analytical expression for it. As a result, Lindley \cite{lindley1980} developed approximate procedures for determining the ratio of two integrals, such as~(\ref{eq6}) and~(\ref{eq12}). Several researchers (\cite{kundu2008}) have used this to derive approximate Bayesian estimates.

      Thus, based on Linux loss functions, we get the Bayesian estimator of $g(\vartheta)=\vartheta$ as 
      \begin{equation}\widehat{\vartheta_{l}} = -\frac{1}{c} \ln E\left(e^{-c \vartheta} \mid Data\right).
      \label{eq13}
      \end{equation}
     By applying Lindley's approximate method in~(\ref{eq13}) for solving ratio of integrals without a closed form, we get the Bayesian estimate of $\vartheta$ as following
     \begin{equation}
\widehat{\vartheta_{l}}=\widehat{\vartheta}_{MLE}-\frac{1}{c}\operatorname{ln}\left[1+\frac{c}{2\mathcal{D}}\left(c\widehat{\vartheta}_{MLE}^2-2a+2\widehat{\vartheta}_{MLE}(b-1)\right)\right].
\label{eq14}
     \end{equation}
\section{Acceptance sampling plans under Type I hybrid censoring scheme}\label{sec3}
In this section, we design ASPs for a lot of units having exponential failure time with PDF given by~(\ref{eq1}). For that, a sample comprising $n$ elements is taken from the lot and is tested  under Type I HCS.

Let $\vartheta_A$ be the acceptable quality level (AQL) and $\vartheta_U$ be the unacceptable quality level (UQL) of a unit item in  the lot. The acceptance and rejection of the lot depends on the probability conditions given by:
\begin{equation}
    \begin{aligned}[b]
      P\left(\text{Lot is rejected} \mid \vartheta=\vartheta_A\right) & \leq \alpha,\\
     P\left(\text{Lot is accepted} \mid \vartheta=\vartheta_U\right) & \leq \beta,
    \end{aligned}
    \label{eq15}
\end{equation}
where $\beta$ is the consumer's risk, $\alpha$ is the producer's risk, and $0<\alpha,\beta<1$.
\subsection{Normal approximation of the Bayesian estimators of $\vartheta$ using delta method}
Here, the Bayesian estimator of the unknown parameter $\vartheta$ is used to design the acceptance sampling plan. In order to proceed, it is necessary to derive the distribution of the Bayesian estimator. Recall that both the Bayesian estimators are derived in terms of $\hat{\vartheta}_{MLE}$, the MLE of $\vartheta$. In \cite{childs2003exact}, Childs et.al derived the exact distribution for the conditional MLE of $\vartheta$, $\widehat{\vartheta}_{MLE}$ for $\mathcal{D} \geq 1$ with PDF as

\begin{multline}
f_{\widehat{\vartheta}_{MLE}}(x)=\left(1-v^n\right)^{-1}\left[ \sum_{\mathcal{D}=1}^{\gamma-1} \sum_{k=0}^\mathcal{D} A_{k, \mathcal{D}} \ q\left(x-\mathcal{T}_{k, \mathcal{D}} ; \frac{\mathcal{D}}{\vartheta}, \mathcal{D}\right)+q\left(x ; \frac{\gamma}{\vartheta}, \gamma\right)\right.\\
+\left.\gamma\left(\begin{array}{l}
n \\
\gamma
\end{array}\right) \sum_{k=1}^\gamma \frac{(-1)^k v^{n-\gamma+k}}{n-\gamma+k}\left(\begin{array}{c}
\gamma-1 \\
k-1
\end{array}\right) q\left(x-T_{k, \gamma} ; \frac{\gamma}{\vartheta}, \gamma\right)\right], \  0<x<n \mathcal{T}
\label{eq16}
\end{multline}

where $$v=\operatorname{exp}\left\{-\frac{\mathcal{T}}{\vartheta}\right\}, \quad T_{k, \mathcal{D}}=\frac{(n-\mathcal{D}+k) \mathcal{T}} {\mathcal{D}},$$ 
$$A_{k, \mathcal{D}}=(-1)^k\left(\begin{array}{l}n \\ \mathcal{D}\end{array}\right)\left(\begin{array}{l}\mathcal{D} \\ k\end{array}\right) u^{n-\mathcal{D}+k},$$ and
$$
q(x ; p,t)= \begin{cases}\frac{p^t}{\Gamma(t)} x^{t-1} e^{-p x}, & x>0 \\ 0, & \text { otherwise. }\end{cases}
$$
From~(\ref{eq16}), the expectation of $\widehat{\vartheta}_{MLE}$ and $E(\widehat{\vartheta}_{MLE}^2)$ are obtained respectively as
\begin{multline}
E(\widehat{\vartheta}_{MLE})=\left(1-v^n\right)^{-1}  {\left[\sum_{\mathcal{D}=1}^{\gamma-1} \sum_{k=0}^\mathcal{D} A_{k, \mathcal{D}}\left(\vartheta+T_{k, \mathcal{D}}\right)+\vartheta\right.} \\
 \left.+\gamma\left(\begin{array}{l}
n \\
\gamma
\end{array}\right) \sum_{k=1}^\gamma \frac{(-1)^k v^{n-\gamma+k}}{n-\gamma+k}\left(\begin{array}{l}
\gamma-1 \\
k-1
\end{array}\right)\left(\vartheta+T_{k, \gamma}\right)\right],
\label{eq17}
\end{multline}
and
\begin{multline}
E\left(\widehat{\vartheta}_{MLE}^2\right)=\left(1-v^n\right)^{-1}\left[ \sum_{\mathcal{D}=1}^{\gamma-1} \sum_{k=0}^\mathcal{D} A_{k, \mathcal{D}}\left\{\frac{\vartheta^2}{\mathcal{D}}(1+\mathcal{D})+2 T_{k, \mathcal{D}} \vartheta+\left(T_{k, \mathcal{D}}\right)^2\right\} \right.\\ +\frac{\vartheta^2}{\gamma}(1+\gamma)
+\gamma\left(\begin{array}{l}
n \\
\gamma
\end{array}\right) \left.\cdot \sum_{k=1}^\gamma \frac{(-1)^k v^{n-\gamma+k}}{n-\gamma+k}\left(\begin{array}{l}
\gamma-1 \\
k-1
\end{array}\right)\left\{\frac{\vartheta^2}{\gamma}(1+\gamma)+2 T_{k, \gamma} \vartheta+\left(T_{k, \gamma}\right)^2\right\}\right].
\label{eq18}
\end{multline}
It is important to note that the variance of  $\widehat{\vartheta}_{MLE}$, denoted as $\sigma^2\left({\widehat{\vartheta}_{MLE}}\right)$, can be obtained by applying equations~(\ref{eq17}) and~(\ref{eq18}).\\

Define the Bayesian estimator of $\vartheta$ using SEL, $\widehat{\vartheta_{s}}=U_1\left(\widehat{\vartheta}_{MLE}\right)$ then by delta method (as in \cite{maheshbajeel2016}), $ \widehat{\vartheta_{s}}$ follows Normal distribution with mean, 
\begin{equation} E\left(\widehat{\vartheta_{s}}\right) = U_1\left(E\left(\widehat{\vartheta}_{MLE}\right)\right) = \frac{\mathcal{D} \ E\left(\hat{\vartheta}_{MLE}\right)+a}{\mathcal{D}+b-1}, 
\label{eq19}
\end{equation} and variance, 
\begin{equation}
\begin{aligned}[b]
V\left(\widehat{\vartheta_{s}}\right) &= \left(U_1^\prime \left(E\left(\widehat{\vartheta}_{MLE}\right)\right)\right)^2\sigma^2\left({\widehat{\vartheta}_{MLE}}\right)\\
  &=\left(\frac{\mathcal{D}}{\mathcal{D}+b-1}\right)^2 \sigma^2\left({\widehat{\vartheta}_{MLE}}\right).
\end{aligned} 
 \label{eq20}
\end{equation} Similarly, from~(\ref{eq14}), we can define the Bayesian estimator of $\vartheta$ using Linex loss function, $\widehat{\vartheta_{l}}$ as $U_2\left(\widehat{\vartheta}_{MLE}\right)$. Then by delta method, $\widehat{\vartheta_{l}}$ follows Normal distribution with mean, 
\begin{equation}
\begin{aligned}[b]
E\left(\widehat{\vartheta_{l}}\right) &= U_2\left(E\left(\widehat{\vartheta}_{MLE}\right)\right)\\
   &= E\left(\widehat{\vartheta}_{MLE}\right)-\frac{1}{c}\operatorname{ln}\left[1+\frac{c}{2\mathcal{D}}\left(cE\left(\widehat{\vartheta}_{MLE}\right)^2-2a+2E\left(\widehat{\vartheta}_{MLE}\right)(b-1)\right)\right], 
   \label{eq21}
\end{aligned} 
\end{equation} and variance,
\begin{equation}
\begin{aligned}[b]
    V\left(\widehat{\vartheta_{l}}\right) &= \left(U_2^\prime \left(E\left(\widehat{\vartheta}_{MLE}\right)\right)\right)^2\sigma^2\left({\widehat{\vartheta}_{MLE}}\right)\\
    &=\left(1-\frac{2c E\left(\widehat{\vartheta}_{MLE}\right) +2b -2}{2\mathcal{D} +c\left(cE\left(\widehat{\vartheta}_{MLE}\right)^2-2a + 2E\left(\widehat{\vartheta}_{MLE}\right) (b-1)\right)}\right)^2 \ \sigma^2\left({\widehat{\vartheta}_{MLE}}\right).
\end{aligned} 
\label{eq22}
\end{equation} 

\subsection{ASP with Bayesian estimator of $\vartheta$ using SEL}
For this ASP, a sample consisting of $n$ items is taken from the lot and examined using Type I hybrid censoring. The probability requirements stated in~(\ref{eq15}) are used to determine whether to accept or reject a lot. Our ASP in this situation is as follows:
\begin{enumerate}
    \item [] Step 1: A random sample of size $n$ is taken from the lot and is tested up to time $\mathcal{T}^*=\operatorname{Min}\{\mathcal{T},X_{\gamma}\}$ where $\mathcal{T}$ is the pre-fixed time and $\gamma < n$ is the pre-fixed number of failures.The number of failures that happened prior to the time $\mathcal{T}^*$ with associated lifetimes provided by~(\ref{eq2}), i.e, either data type I or type II, is denoted by $\mathcal{D}$. Also, fix the values of AQL, UQL, producer's risk ($\alpha$), and consumer's risk ($\beta$).
    \item [] Step 2: Calculate the Bayesian estimate of $\vartheta$ using SEL ($ \widehat{\vartheta_{s}}$), which is given by 
     \begin{equation*}
       \widehat{\vartheta_{s}} = \frac{\mathcal{D} \ \hat{\vartheta}_{MLE}+a}{\mathcal{D}+b-1} .
      \end{equation*}
      \item[] Step 3: Continue the testing process, by repeating Steps 1, and  2 if $t_1 \leq  \widehat{\vartheta_{s}} < t_2$. If not, go to Step 4.
      \item[] Step 4: Accept the lot if $ \widehat{\vartheta_{s}} \geq t_2$ and reject the lot if $ \widehat{\vartheta_{s}}< t_1.$
    \end{enumerate}

The key objective is to minimize the total testing cost subject to the constrains given in~(\ref{eq15}) in order to obtain the optimal values of $n, t_1$ and $t_2$. The total testing cost is the product of the testing cost of an item for unit time with total testing time. Based on the sampling plan described above, the total testing time is the product of the time taken to reach a decision for each sample, $\widehat{\vartheta_{s}}$, by the number of samples tested. But, $\widehat{\vartheta_{s}}$ and number of samples are random variables. So, here the expected testing cost is obtained.

Consider $p_c$, $p_a$, and $p_r$ as the probabilities of continuation of testing, accepting the lot and rejecting the lot respectively. Then, the corresponding long run probabilities of acceptance and rejection are given by $P_a = \frac{p_a}{1-p_c}$ and $P_r = \frac{p_r}{1-p_c}$, and expected number of items failed is given by $\frac{1}{1-p_c}$ (see \cite{sherman1965}). Thus, the expected testing cost for the ASP is obtained as $$\text{ETC} = \frac{C \ E\left(\widehat{\vartheta_{s}}\right)}{1-p_c},$$
where $C$ is the testing cost of an item for unit time.

From equations~(\ref{eq19}) and~(\ref{eq20}), we get that, $\widehat{\vartheta_{s}}$ follows normal distribution with mean, $$E\left(\widehat{\vartheta_{s}}\right)=\frac{\mathcal{D} \ E\left(\hat{\vartheta}_{MLE}\right)+a}{\mathcal{D}+b-1},$$ and variance, $$V\left(\widehat{\vartheta_{s}}\right)=\left(\frac{\mathcal{D}}{\mathcal{D}+b-1}\right)^2 \sigma^2\left({\widehat{\vartheta}_{MLE}}\right).$$ Then, the required probabilities are given by,
\begin{equation}
    p_a = P\left(\widehat{\vartheta_{s}} \geq t_2\right) = P\left[Z \geq \frac{t_2 - E\left(\widehat{\vartheta_{s}}\right)}{\sqrt{V\left(\widehat{\vartheta_{s}}\right)}}\right],
    \label{eq23}
\end{equation}
\begin{equation}
    p_r = P\left(\widehat{\vartheta_{s}} < t_1\right) = P\left[Z < \frac{t_1 - E\left(\widehat{\vartheta_{s}}\right)}{\sqrt{V\left(\widehat{\vartheta_{s}}\right)}}\right],
    \label{eq24}
\end{equation}
\begin{equation}
    p_c = P\left(t_1 \leq \widehat{\vartheta_{s}} < t_2\right) = P\left[\frac{t_1 - E\left(\widehat{\vartheta_{s}}\right)}{\sqrt{V\left(\widehat{\vartheta_{s}}\right)}} \leq Z < \frac{t_2 - E\left(\widehat{\vartheta_{s}}\right)}{\sqrt{V\left(\widehat{\vartheta_{s}}\right)}}\right].
    \label{eq25}
\end{equation}

    As the expected testing cost is a function of the unknown parameter $\vartheta$, here we consider the expected testing cost at $\vartheta_A$ as the objective of the optimisation problem. Thus, the required optimisation problem, $P_1\left(n,t_1,t_2\right)$, is obtained by minimising the expected testing cost at $\vartheta_A$, subject to the probability conditions in~(\ref{eq15}). 
    \begin{equation} 
\begin{aligned}[b]
P_1(n, t_1, t_2):\qquad &\underset{\left(n, t_{1}, t_{2}\right)}{\operatorname{Min}} \left(\frac{C \ E\left(\widehat{\vartheta_{s}}\right)}{1-p_{c}}\right) \text{at} \ \vartheta_A \\
\text {subject to} \qquad &\left( P_r \mid \vartheta = \vartheta_A\right) \leq \alpha, \\
&\left( P_a \mid \vartheta = \vartheta_U \right) \leq \beta, 
\end{aligned}
\label{eq26}
\end{equation}
where $t_{1}, t_{2}>0$ and $t_{2}>t_{1}$. That is,
\begin{equation}
\begin{aligned}[b]
\quad&\underset{\left(n, t_{1}, t_{2}\right)}{\operatorname{Min}} \left(\frac{C \ E\left(\widehat{\vartheta_{s}}\right)}{1-p_{c}}\right) \text{at} \ \vartheta_A \\
\text {subject to} \qquad &\left(\frac{p_{r}}{1-p_{c}} \mid \vartheta = \vartheta_A\right) \leq \alpha, \\
&\left(\frac{p_{a}}{1-p_{c}} \mid \vartheta = \vartheta_U \right) \leq \beta,
\end{aligned}
\label{eq27}
\end{equation}
where $t_{1}, t_{2}>0$ and $t_{2}>t_{1}$. Thus, by substituting the expressions of $p_{a}, p_{r}$ and $p_{c}$ from equations~(\ref{eq23}),~(\ref{eq24}) and~(\ref{eq25}), the non-linear optimization problem, $P_1(n,t_1, t_2)$, becomes:
  \begin{equation} 
\begin{aligned}[b]
P_1(n, t_1, t_2):\qquad &\underset{\left(n, t_{1}, t_{2}\right)}{\operatorname{Min}} \left(\frac{C \ E\left(\widehat{\vartheta_{s}}\right)}{1-P\left[\frac{t_1 - E\left(\widehat{\vartheta_{s}}\right)}{\sqrt{V\left(\widehat{\vartheta_{s}}\right)}} \leq Z < \frac{t_2 - E\left(\widehat{\vartheta_{s}}\right)}{\sqrt{V\left(\widehat{\vartheta_{s}}\right)}}\right]}\right) \text{at} \ \vartheta_A \\
\text {subject to} \qquad &\left( \frac{P\left[Z < \frac{t_1 - E\left(\widehat{\vartheta_{s}}\right)}{\sqrt{V\left(\widehat{\vartheta_{s}}\right)}}\right]}{1-P\left[\frac{t_1 - E\left(\widehat{\vartheta_{s}}\right)}{\sqrt{V\left(\widehat{\vartheta_{s}}\right)}} \leq Z < \frac{t_2 - E\left(\widehat{\vartheta_{s}}\right)}{\sqrt{V\left(\widehat{\vartheta_{s}}\right)}}\right]} \mid \vartheta = \vartheta_A\right) \leq \alpha, \\
&\left( \frac{ P\left[Z \geq \frac{t_2 - E\left(\widehat{\vartheta_{s}}\right)}{\sqrt{V\left(\widehat{\vartheta_{s}}\right)}}\right]}{1-P\left[\frac{t_1 - E\left(\widehat{\vartheta_{s}}\right)}{\sqrt{V\left(\widehat{\vartheta_{s}}\right)}} \leq Z < \frac{t_2 - E\left(\widehat{\vartheta_{s}}\right)}{\sqrt{V\left(\widehat{\vartheta_{s}}\right)}}\right]}\mid \vartheta = \vartheta_U \right) \leq \beta, 
\end{aligned}
\label{eq28}
\end{equation}
where $t_{1}, t_{2}>0$ and $t_{2}>t_{1}$.

For $n, t_1,$ and $ t_2$, this non-linear optimisation problem $P_1$ can be solved by  the following steps:
\begin{enumerate}
    \item  Compute the minimum value of $\gamma$ which satisfies the constrains of the problem $P_1$.
    \item Utilising the obtained value of $\gamma$, from Step 1, solve $P_1$ for the optimal values of $n$, $t_1$, and $t_2$.
\end{enumerate}
The optimal values of $n$, $t_1$, and $t_2$ are obtained by solving the above non-linear optimization problem using genetic algorithm solver in MATLAB and are tabulated in Table~\ref{Tab:SEL}.
\subsection{ASP with Bayesian estimator of $\vartheta$ using Linex loss function}
In this section, the Bayesian estimate of $\vartheta$ using Linex loss function is applied in the ASP. For that, a sample of size $n$ is taken from the lot. The testing is terminated at time $\mathcal{T}^*$ and the observed lifetime data will be as given in~(\ref{eq2}). Here we define an ASP as follows: 
\begin{enumerate}
    \item [] Step 1: Take a random sample of size $n$  from the lot. Fix the time $\mathcal{T}$, failure number, $\gamma < n$, AQL, UQL, producer's risk ($\alpha$), and consumer's risk ($\beta$). Testing is stopped at time, $\mathcal{T}^*=\operatorname{Min}\{\mathcal{T},X_{\gamma}\}$.
    \item [] Step 2: Calculate the Bayesian estimate of $\vartheta$ using Linex loss function ($ \widehat{\vartheta_{l}}$),  given by 
     \begin{equation*}
      \widehat{\vartheta_{l}}=\widehat{\vartheta}_{MLE}-\frac{1}{c}\operatorname{ln}\left[1+\frac{c}{2\mathcal{D}}\left(c\widehat{\vartheta}_{MLE}^2-2a+2\widehat{\vartheta}_{MLE}(b-1)\right)\right].  
     \end{equation*}
      \item[] Step 3: Continue the test, if $t_1 \leq  \widehat{\vartheta_{l}} < t_2$ and repeat Steps 1 and 2 .\\ Else, go to Step 4.
      \item[] Step 4: Accept the lot, if $ \widehat{\vartheta_{l}} \geq t_2$ and reject the lot, if $ \widehat{\vartheta_{l}}< t_1.$
    \end{enumerate}
 The main problem here is determining the optimal values of $n$, $t_1$, and $t_2$ while minimising the testing cost under probability constraints as mentioned in~(\ref{eq15}). Total testing cost is a product of testing cost of a unit item in unit time and total testing time. Total testing time is obtained by multiplying the overall number of samples tested with decision making time while testing a sample, given by  $\widehat{\vartheta}_l$. But they are both random variables. Thus, the total testing time becomes a random quantity  with expectation, $\frac{1}{1-p_c} \ E\left(\widehat{\vartheta_{l}}\right)$. As a result, here we use the expected testing cost as the objective function and it is given by
 $$\text{ETC} = \frac{C \ E\left(\widehat{\vartheta_{l}}\right)}{1-p_c},$$ where $C$ is the testing cost of an item for unit time.

 From equations~(\ref{eq21}) and~(\ref{eq22}), we get that, $\widehat{\vartheta_{l}}$ follows Normal distribution with mean, $$E\left(\widehat{\vartheta_{l}}\right) =  E\left(\widehat{\vartheta}_{MLE}\right)-\frac{1}{c}\operatorname{ln}\left[1+\frac{c}{2\mathcal{D}}\left(cE\left(\widehat{\vartheta}_{MLE}\right)^2-2a+2E\left(\widehat{\vartheta}_{MLE}\right)(b-1)\right)\right],$$ and variance, $$V\left(\widehat{\vartheta_{l}}\right) =\left(1-\frac{2c E\left(\widehat{\vartheta}_{MLE}\right) +2b -2}{2\mathcal{D} +c\left(cE\left(\widehat{\vartheta}_{MLE}\right)^2-2a + 2E\left(\widehat{\vartheta}_{MLE}\right) (b-1)\right)}\right)^2 \ \sigma^2\left({\widehat{\vartheta}_{MLE}}\right). $$ Using this information, we can evaluate the required probabilities as
 \begin{equation}
    p_a = P\left(\widehat{\vartheta_{l}} \geq t_2\right) = P\left[Z \geq \frac{t_2 - E\left(\widehat{\vartheta_{l}}\right)}{\sqrt{V\left(\widehat{\vartheta_{l}}\right)}}\right],
    \label{eq29}
\end{equation}
\begin{equation}
    p_r = P\left(\widehat{\vartheta_{l}} < t_1\right) = P\left[Z < \frac{t_1 - E\left(\widehat{\vartheta_{l}}\right)}{\sqrt{V\left(\widehat{\vartheta_{l}}\right)}}\right],
    \label{eq30}
\end{equation}
\begin{equation}
    p_c = P\left(t_1 \leq \widehat{\vartheta_{l}} < t_2\right) = P\left[\frac{t_1 - E\left(\widehat{\vartheta_{l}}\right)}{\sqrt{V\left(\widehat{\vartheta_{l}}\right)}} \leq Z < \frac{t_2 - E\left(\widehat{\vartheta_{l}}\right)}{\sqrt{V\left(\widehat{\vartheta_{l}}\right)}}\right].
    \label{eq31}
\end{equation}
 
Here, we evaluate the objective function of the optimisation problem at $\vartheta_A$ since the overall testing cost includes the unknown parameter $\vartheta$. Thus, the expected testing cost at $\vartheta_A$ is minimised by satisfying the probability conditions given in~(\ref{eq15}). Hence, the required non-linear optimisation problem is formulated as
\begin{equation} 
\begin{aligned}[b]
P_2(n, t_1, t_2):\qquad &\underset{\left(n, t_{1}, t_{2}\right)}{\operatorname{Min}} \left(\frac{C \ E\left(\widehat{\vartheta_{l}}\right)}{1-p_{c}}\right) \text{at} \ \vartheta_A\\
\text {subject to} \qquad &\left( P_r \mid \vartheta = \vartheta_A\right) \leq \alpha, \\
&\left( P_a \mid \vartheta = \vartheta_U \right) \leq \beta, 
\end{aligned}
\label{eq32}
\end{equation}
where $t_{1}, t_{2}>0$ and $t_{2}>t_{1}$. That is,
\begin{equation}
\begin{aligned}[b]
\quad&\underset{\left(n, t_{1}, t_{2}\right)}{\operatorname{Min}} \left(\frac{C \ E\left(\widehat{\vartheta_{l}}\right)}{1-p_{c}}\right) \text{at} \ \vartheta_A \\
\text {subject to} \qquad &\left(\frac{p_{r}}{1-p_{c}} \mid \vartheta = \vartheta_A\right) \leq \alpha, \\
&\left(\frac{p_{a}}{1-p_{c}} \mid \vartheta = \vartheta_U \right) \leq \beta,
\end{aligned}
\label{eq33}
\end{equation}
where $t_{1}, t_{2}>0$ and $t_{2}>t_{1}$. The non-linear optimisation problem, $P_2(n,t_1, t_2)$, is modified by substituting the expressions for $p_a$, $p_r$, and $p_c$ from equations~(\ref{eq29}),~(\ref{eq30}), and~(\ref{eq31}) as 
  \begin{equation} 
\begin{aligned}[b]
P_2(n, t_1, t_2):\qquad &\underset{\left(n, t_{1}, t_{2}\right)}{\operatorname{Min}} \left(\frac{C \ E\left(\widehat{\vartheta_{l}}\right)}{1-P\left[\frac{t_1 - E\left(\widehat{\vartheta_{l}}\right)}{\sqrt{V\left(\widehat{\vartheta_{l}}\right)}} \leq Z < \frac{t_2 - E\left(\widehat{\vartheta_{l}}\right)}{\sqrt{V\left(\widehat{\vartheta_{l}}\right)}}\right]}\right) \text{at} \ \vartheta_A \\
\text {subject to} \qquad &\left( \frac{P\left[Z < \frac{t_1 - E\left(\widehat{\vartheta_{l}}\right)}{\sqrt{V\left(\widehat{\vartheta_{l}}\right)}}\right]}{1-P\left[\frac{t_1 - E\left(\widehat{\vartheta_{l}}\right)}{\sqrt{V\left(\widehat{\vartheta_{l}}\right)}} \leq Z < \frac{t_2 - E\left(\widehat{\vartheta_{l}}\right)}{\sqrt{V\left(\widehat{\vartheta_{l}}\right)}}\right]} \mid \vartheta = \vartheta_A\right) \leq \alpha, \\
&\left( \frac{ P\left[Z \geq \frac{t_2 - E\left(\widehat{\vartheta_{l}}\right)}{\sqrt{V\left(\widehat{\vartheta_{l}}\right)}}\right]}{1-P\left[\frac{t_1 - E\left(\widehat{\vartheta_{l}}\right)}{\sqrt{V\left(\widehat{\vartheta_{l}}\right)}} \leq Z < \frac{t_2 - E\left(\widehat{\vartheta_{l}}\right)}{\sqrt{V\left(\widehat{\vartheta_{l}}\right)}}\right]}\mid \vartheta = \vartheta_U \right) \leq \beta, 
\end{aligned}
\label{eq34}
\end{equation}
where $t_{1}, t_{2}>0$ and $t_{2}>t_{1}$.

The following steps are used to resolve the above mentioned nonlinear optimisation problem $P_2$ for ($n, t_1, t_2$):
\begin{enumerate}
    \item First, determine the least value of $\gamma$ subjected to the constrains of the problem $P_2$.
    \item Then solve $P_2$, using the value of $\gamma$ obtained in Step 1.
\end{enumerate}
This non-linear optimisation problem is solved using genetic algorithm solver in MATLAB and the obtained optimal values of $n$, $t_1$, and $t_2$ are tabulated in Table~\ref{Tab:LL1} and Table~\ref{Tab:LL2} for $c=0.5$ and $c=-0.5$, respectively.
\section{Computational results, comparisons and real data case study}\label{sec4}
Numerical computations of all the above sections are discussed in this section and are displayed in Table~\ref{Tab:SEL} to \ref{Tab:LL2}.  The hyper parameters, $a=1.25$ and $b=2.5$ are used to compute the Bayesian estimator throughout this work since they give more approximate Bayesian estimate for $\vartheta$ (see \cite{lin2008exact} and \cite{yeh1994bayesian}). ETC of the ASP under Type I hybrid censoring using Bayesian estimate of $\vartheta$ subject to Linex loss function for both $c=0.5$ and $c=-0.5$ is computed and is shown in Table~\ref{Tab:LL1}  and Table~\ref{Tab:LL2}, respectively. From that, we can observe that, the ETC is less for $c=0.5$ than that of $c=-0.5$. For example, with $C=1$, $\vartheta_A=200$, $\vartheta_U=100$, $\mathcal{T}=100$, $\alpha=0.05$ and $\beta=0.05$, gives ETC of 127.2029 for $c=0.5$ and ETC is increased to 228.1767 for $c=-0.5$. A comparison of testing costs obtained from ASPs under Type I hybrid censoring using a Bayesian estimate of the parameter $\vartheta$ using the SEL and Linex loss functions is also done using the results in Tables~\ref{Tab:SEL}, \ref{Tab:LL1} and \ref{Tab:LL2}. Consider the following example for the set of parameters: $C=1$, $\vartheta_A=500$, $\vartheta_U=200$, $\mathcal{T}=50$, $\alpha=0.05$ and $\beta=0.05$, the ETC obtained for ASP with SEL and Linex loss functions are 475.5810, 557.9282 ($c=0.5$) and 627.7466 ($c=-0.5$), respectively. These comparisons shows that the ASP has the lowest test cost when Bayesian estimate of $\vartheta$ is computed using SEL function. 

Table~\ref{Tab:comparison} illustrates the comparison of ETC computed for the ASP under Type I hybrid censoring with Bayesian estimate of $\vartheta$ using SEL function with the ASP under Type I censoring scheme designed using MLE of $\vartheta$ in \cite{kumar2016design}. It can be seen that, for parameters, $C=1$, $\vartheta_A=200$, $\vartheta_U=100$, $\mathcal{T}=100$, $\alpha=0.05$ and $\beta=0.25$, the ETC of ASP using SEL is 144.1837 and 156.53 for the ASP in \cite{kumar2016design}. Thus, one can observe that cost of our plan is less than the one described in \cite{kumar2016design}. Also, the proposed acceptance sampling plans are applied to a real life data. 
\begin{table}[ht]
\centering
{\tiny
\begin{tabular}{llllllllll}
\hline
&&&&&&&&&\\
$\vartheta_A$ & $\vartheta_U$ &$\mathcal{T}$ &$\alpha$&$\beta$&$\gamma$& $t_1$ &$t_2$&$n$&ETC\\
&&&&&&&&&\\
\hline 
&&&&&&&&&\\
200&100&100&0.05&0.05&26&162.3926&162.3957&31&123.4077\\
&&&&&&&&&\\
&&&0.01&0.05&20&74.7316&74.7322&23&90.7417\\
&&&&&&&&&\\
&&&0.01&0.01&25&146.3421&146.3394&26&125.6045\\
&&&&&&&&&\\
500&200&50&0.05&0.05&21&313.5638&313.5638&26&475.5810\\
&&&&&&&&&\\
&&&0.01&0.05&22&310.5940&310.5948&27&480.9310\\
&&&&&&&&&\\
&&&0.01&0.01&30&225.7734&225.7735&32&563.9248\\
&&&&&&&&&\\
500&200&100&0.05&0.05&26&313.3180&313.3188&37&509.9392\\
&&&&&&&&&\\
&&&0.01&0.05&14&311.1179&311.1179&19&503.4323\\
&&&&&&&&&\\
&&&0.01&0.01&24&390.4518&390.4526&38&566.6305\\
&&&&&&&&&\\
3000&1500&1000&0.05&0.05&16&2005.3416&2005.3425&62&475.5810\\
&&&&&&&&&\\
&&&0.01&0.05&10&2000.6681&2000.6682&92&480.9310\\
&&&&&&&&&\\
&&&0.01&0.01&15&2022.5679&2022.568&72&2684.0441\\
&&&&&&&&&\\
 \hline
\end{tabular}}
\caption{ASP under Type I HCS using SEL function for C=1.}
\label{Tab:SEL}
\end{table}

\begin{table}[ht]
\centering
{\tiny
\begin{tabular}{llllllllll}
\hline
&&&&&&&&&\\
$\vartheta_A$ & $\vartheta_U$ &$\mathcal{T}$ &$\alpha$&$\beta$&$\gamma$& $t_1$ &$t_2$&$n$&ETC\\
&&&&&&&&&\\
\hline 
&&&&&&&&&\\
200&100&100&0.05&0.05&21&63.1133&63.1134&35&127.2029\\
&&&&&&&&&\\
&&&0.01&0.05&19&56.3674&56.3675&50&118.7482\\
&&&&&&&&&\\
&&&0.01&0.01&28&140.9731&140.9733&30&143.1099\\
&&&&&&&&&\\
500&200&50&0.05&0.05&21&367.4209&367.421&34&557.9282\\
&&&&&&&&&\\
&&&0.01&0.05&26&361.9554&361.9555&32&556.9578\\
&&&&&&&&&\\
&&&0.01&0.01&28&371.9294&371.9295&36&595.2252\\
&&&&&&&&&\\
500&200&100&0.05&0.05&22&367.4918&367.4918&39&521.9267\\
&&&&&&&&&\\
&&&0.01&0.05&27&354.8226&354.8229&37&513.5307\\
&&&&&&&&&\\
&&&0.01&0.01&25&364.1113&364.1115&36&524.21\\
&&&&&&&&&\\
3000&1500&1000&0.05&0.05&15&2473.6173&2473.6174&24&3379.9555\\
&&&&&&&&&\\
&&&0.01&0.05&11&2437.8132&2437.8133&33&3161.2578\\
&&&&&&&&&\\
&&&0.01&0.01&6&2428.6943&2428.6944&14&3705.1855\\
&&&&&&&&&\\
 \hline
\end{tabular}}
\caption{ASP under Type I HCS using Linex Loss for $c=0.5$ and C=1.}
\label{Tab:LL1}
\end{table}

\begin{table}[ht]
\centering
{\tiny
\begin{tabular}{llllllllll}
\hline
&&&&&&&&&\\
$\vartheta_A$ & $\vartheta_U$ &$\mathcal{T}$ &$\alpha$&$\beta$&$\gamma$& $t_1$ &$t_2$&$n$&ETC\\
&&&&&&&&&\\
\hline 
&&&&&&&&&\\
200&100&100&0.05&0.05&28&186.2786&186.2793&34&228.1767\\
&&&&&&&&&\\
&&&0.01&0.05&26&178.3335&178.3342&29&232.5816\\
&&&&&&&&&\\
&&&0.01&0.01&23&177.5453&177.5455&45&220.3530\\
&&&&&&&&&\\
500&200&50&0.05&0.05&12&400.025&400.0248&25&627.7466\\
&&&&&&&&&\\
&&&0.01&0.05&15&396.9158&396.9158&26&635.9504\\
&&&&&&&&&\\
&&&0.01&0.01&14&405.41&405.4198&26&636.7619\\
&&&&&&&&&\\
500&200&100&0.05&0.05&18&423.5114&423.5116&43&600.3209\\
&&&&&&&&&\\
&&&0.01&0.05&22&405.5332&405.5335&38&620.4932\\
&&&&&&&&&\\
&&&0.01&0.01&24&411.5643&411.5644&37&623.4914\\
&&&&&&&&&\\
3000&1500&1000&0.05&0.05&19&3783.5928&3783.5929&29&3449.2903\\
&&&&&&&&&\\
&&&0.01&0.05&11&2507.4709&2507.471&30&3411.8618\\
&&&&&&&&&\\
&&&0.01&0.01&9&2557.26&2557.2699&27&3450.2116\\
&&&&&&&&&\\
 \hline
\end{tabular}}
\caption{ASP under Type I HCS using Linex Loss for $c=-0.5$ and C=1.}
\label{Tab:LL2}
\end{table}
\begin{table}[ht]
\centering
\begin{tabular}{llllllll}
\hline
$\vartheta_A$ & $\vartheta_U$ &$\mathcal{T}$&$\alpha$&$\beta$&$A$&$B$&$C$\\
\hline
200&100&100&0.05&0.25&123.4077&151.0479&156.53\\
3000&1500&1500&0.05&0.1&1423.5748&2256.9985&2290.35\\
 \hline
\end{tabular}
\caption{ Comparison of ETCs for ASP using two different loss functions and testing cost under Type I Censoring.  The ETCs are calculated using the SEL function in column A and Linex loss with $c=0.5$ in column B. Column C presents the expected testing cost obtained for ASP under Type I censoring, as explained in Kumar et al.'s study in \cite{kumar2016design}. $C=1$}
\label{Tab:comparison}
\end{table}
\subsection{Case Study}
In this section, the acceptance sampling plans presented in this work are demonstrated using a real-world example from \cite{alizadeh2015testing}. The failure rates for 36 appliances that underwent an automated life test are listed below; the lifetime shown here refers to the number of cycles of use that the appliances can withstand before failing. 
\begin{multline*}
\begin{aligned}
    &11, 35, 49, 170, 329, 381, 708, 958, 1062, 1167, 1594, 1925, 1990, 2223, 2327,
    2400, 2451, 2471,2551, 2565, 2568,\\
    &2694, 2702, 2761, 2831, 3034, 3059, 3112,
3214, 3478, 3504, 4329, 6367, 6976, 7846, 13403.
    \end{aligned}
\end{multline*}
The above data follows exponential distribution. In order to apply the proposed sampling plans in Sections 3 and 4, we take our plan parameters as $\vartheta_A=3000$, $\vartheta_U=600$, $\mathcal{T}=2000$, $\alpha=0.1$ and $\beta=0.2$ as in \cite{kumar2016design}. Now, we can solve the nonlinear optimisation problems $P_1$ and $P_2$ to obtain the optimal values of plan parameters $t_1$ and $t_2$. Note that, $t_1$ and $t_2$ are real numbers but for the application of our acceptance sampling plans we need them to be integers. So we first solve $P_1$ and $P_2$ , if we get integer values for $t_1$ and $t_2$, then we can look for the acceptance and rejection of lot with appropriate conditions. If $t_1$ and $t_2$ are non-integers then we restrict them to be integers. Hence our optimisation problems $P_1$ and $P_2$ becomes non-linear integer programming problem with integer solutions. 

For the above mentioned data, the different sampling plans proposed in our work yield the following outcomes:

\begin{enumerate}
    \item ASP under Type I HCS using SEL function gives the parameters as $\gamma=9$, $t_1 = 2064$, $t_2=2065$, and $n=31$. Next, we consider the first $31$ samples from the above data. The samples are tested upto time $\mathcal{T}^*=\operatorname{Min}\{\mathcal{T},X_{9}\} =1062$ since $X_9=1062$. Now we compute the Bayesian estimate of $\vartheta$ with SEL function as $\widehat{\vartheta_{s}}=2577.9286$ and it is greater than $t_2=2065$. So, we can accept the lot with ETC = 2405.
    \item Based on ASP under Type I hybrid using Linex loss function for $c=0.5$, we get the optimal values for test parameters as $\gamma=11$, $t_1 = 2156$, $t_2=2157$, and $n=27$. Then, the first $27$ samples are taken from the data and the testing is stopped at  $\mathcal{T}^*=\operatorname{Min}\{\mathcal{T},X_{11}\} =1594$ as $X_{11}=1594$. Next, we calculate the Bayesian estimate of $\vartheta$ using Linex loss function and get $\widehat{\vartheta_{l}}=2883.2339$. Hence, we accept the lot as $\widehat{\vartheta_{l}}>t_2$ and the ETC = 2909.
\end{enumerate}
\section{Conclusion}\label{sec5}
In conclusion, a novel approach for Type I hybrid censoring in acceptance sampling plans is proposed here by utilizing Bayesian inference. Unlike conventional methods that rely on maximum likelihood estimation (MLE) of $\vartheta$ and minimizing Bayes risk, this study considered the Bayesian estimators of $\vartheta$, for designing the ASP and taking into account the expected testing cost. The squared error loss and linex loss were employed to compute the Bayesian estimators of $\vartheta$. 

The aim of all ASPs was to compute the optimal values of $n$, $t_1$ and $t_2$ by solving nonlinear optimisation problems of minimising ETC. A comparative study of ETC involved in ASPs using SEL and Linex loss functions is illustrated here using Table~\ref{Tab:SEL} to \ref{Tab:LL2} for different values of AQL ($\vartheta_A$), UQL ($\vartheta_U$), producer's risk ($\alpha$), consumer's risk ($\beta$), $\mathcal{T}$ and $\gamma$.
The lower expected testing cost for SEL compared to Linex loss functions signifies that SEL provides a more cost-effective approach for our plan. This insight can have substantial implications for resource allocation, budget optimization, and overall project efficiency.
ETC obtained for ASP using Linex loss function is computed for both $c=0.5$ and $c=-0.5$ and is tabulated in Table~\ref{Tab:LL1} and Table~\ref{Tab:LL2}, respectively. From that, we get the ETC for $c=-0.5$ is more than ETC obtained for $c=0.5$. This is mainly because the Bayesian estimator of $\vartheta$ is overestimated for $c=-0.5$ using Linex loss function (see \cite{peng2013bayesian}). 

Table~\ref{Tab:comparison}  shows the comparison of expected costs for our plan using Bayesian inference and the  ASP under Type I censoring in \cite{kumar2016design} which uses MLE of $\vartheta$. The comparison results has revealed the importance of utilizing Bayesian inference to minimize costs as the ETC involved in ASP under Type I censoring in \cite{kumar2016design} is more than the ASPs discussed in this work. Thus, by adopting Bayesian inference and incorporating SEL as the loss function, we can better manage costs and improve the overall performance of our plan. Moreover, a real life case study is conducted to explain the practical application of the various ASPs discussed here. In terms of future research possibilities, the application of Bayesian inference can be extended to the design of ASPs, incorporating accelerated life testing, which helps to reduce the testing time.

\bibliographystyle{is-plain}
  \bibliography{references}
\end{document}